\documentclass[12pt]{article}
\usepackage{amsmath}
\usepackage{amssymb}
\usepackage{graphicx}

\newtheorem{theorem}{Theorem}
\newtheorem{lemma}{Lemma}

\title{Exponential Stabilization of Nonholonomic Systems by means of Oscillating Controls\thanks{This work is supported in part by the Alexander von Humboldt Foundation and by the agreement on scientific cooperation between the National Academy of Sciences of Ukraine
and the Polish Academy of Sciences.}}

\date{\small Max Planck Institute for Dynamics of Complex Technical Systems, \\ Magdeburg, Germany}

\author{Alexander Zuyev\thanks{On leave from the the Institute of Applied Mathematics and Mechanics, National Academy of Sciences of Ukraine, Donetsk.
({zuyev@mpi-magdeburg.mpg.de}, {alexander.zuyev@gmail.com}).}}

\begin{document}
\maketitle

\begin{abstract}
This paper is devoted to the stabilization problem for nonlinear driftless control systems by means of a time-varying feedback control.
It is assumed that the vector fields of the system together with their first order Lie brackets span the whole tangent space at the equilibrium. A family of trigonometric open-loop controls is constructed to approximate the gradient flow associated with a Lyapunov function.
These controls are applied for the derivation of a time-varying feedback law under the sampling strategy.
By using Lyapunov's direct method, we prove that the controller proposed ensures exponential stability of the equilibrium.
As an example, this control design procedure is applied to stabilize the Brockett integrator.
\end{abstract}

\section{Introduction}

Consider a control system
\begin{equation}
\dot x = \sum_{i=1}^m u_i f_i (x)\equiv f(x,u),
\label{Sigma}
\end{equation}
where $x=(x_1,x_2,...,x_n)'\in D\subset \mathbb R^n$ is the state and $u=(u_1,u_2,...,u_m)'\in\mathbb R^m$ is the control. The domain $D$ contains the trivial equilibrium point $x=0$.
We treat all vectors as columns and denote the transpose with a prime. The vector fields $f_i(x)$
are assumed to be mappings of class $C^2$ from $D$ to $\mathbb R^n$.

It is a well-known fact due to R.W.~Brockett~\cite{Br1983} that system~\eqref{Sigma} is not stabilizable by a smooth feedback law $u=k(x)$ such that $k(0)=0$, provided that $m<n$ and $f_1(0)$, $f_2(0)$, ...., $f_m(0)$ are linearly independent vectors. Note that Brockett's condition remains necessary for the stabilizability in a class of discontinuous feedback laws provided that the solutions of the closed-loop system are defined in the sense of A.F.~Filippov~\cite{Ryan}.
To overcome this obstruction, two main strategies can be used for the stabilization of general controllable systems.
The first strategy is based on the use of a time-varying continuous feedback law $u=k(t,x)$ to stabilize the origin of a small-time locally controllable system~\cite{Coron}. In the other strategy, the equilibrium of an asymptotically controllable system can be stabilized by means of a discontinuous feedback law $u=k(x)$, provided that the solutions (``$\pi$-trajectories'') are defined in the sense of sampling~\cite{CLSS}.

An approach for the practical stabilization of nonholomomic systems based on transverse functions is proposed by P.~Morin and C.~Samson~\cite{MorS}. A survey of feedback design techniques is presented in the book by J.-M.~Coron~\cite{Coron}. Despite the rich literature in this area and to the best of our knowledge, there is no universal procedure available for the stabilizing control design for an arbitrary nonlinear system of form~\eqref{Sigma}.

The paper~\cite{CdWS} is devoted to the control design for a kinematic cart model with two inputs. A coordinate transformation from the three-dimensional state space to a two-dimensional manifold (parameterized by the arc length and the orientation error) plays a crucial role in the analysis. Based on this representation, a discontinuous  feedback law is proposed such that any solution of the closed-loop system exponentially converges to an equilibrium point. The orientation angle is defined modulo $2\pi$ in such equilibria.

Applications of sinusoidal controls to the steering problem for systems of form~\eqref{Sigma} are considered in the paper~\cite{MS}.
A combination of constant controls and sinusoids at integrally related frequencies is used to steer the first-order canonical system to an arbitrary configuration. Some modifications of this algorithm are presented for chained systems. An overview of algorithms for the motion planning of nonholonomic systems is presented in the book~\cite{Jean}.

In the paper~\cite{Ailon}, the controllability and trajectory tracking problems are considered for a kinematic car model with nonholonomic constraints. A result on the solvability of the motion planning problem is established for such a model by using trigonometric controls. The error dynamics in a neighborhood of the reference trajectory is studied to solve the tracking problem. It is shown that the error dynamics is stabilizable by using a quadratic Lyapunov function. The controller design scheme proposed is illustrated by examples of a state-to-state control and tracking a circle with time scheduling at selected points.

The stabilization problem for a nonholonomic system in power form with bounded inputs is considered in the paper~\cite{AM}. The receding-horizon principle is used to solve an open-loop optimization problem and to derive a sampling control. It is proved that the family of controls obtained can be used to stabilize the destination state in finite time with any chosen precision. The numerical implementation of this algorithm is shown for a five dimensional system.

The paper~\cite{WTSML} is devoted to the stabilization problem of nonholonomic systems about a feasible trajectory, instead of a point. For such kind of problem, a time-varying feedback law is obtained by using the linearization around a feasible trajectory. The Heisenberg system and a mobile robot model are considered as examples for stabilizing a straight line trajectory in the three-dimensional space. This approach is shown to be applicable for the trajectory stabilization of a front wheel drive car.

Assume that $m<n$ and that $f_1(x)$, $f_2(x)$,..., $f_m(x)$ together with a fixed set of the first order Lie brackets span the whole tangent space for system~\eqref{Sigma}, i.e.
\begin{equation}
{\rm span}\,\{f_i(x), [f_j,f_l](x)\, | \,i=1,2,..., m,\; (j,l)\in S \} = \mathbb R^n,
\label{rank}
\end{equation}
for each $x\in D$,
where $S\subseteq \{1,2,...,m\}^2$,
$$
[f_j,f_l](x) = \frac{\partial f_l(x)}{\partial x}f_j(x) - \frac{\partial f_j(x)}{\partial x}f_l(x)
$$
and $\frac{\partial f_j(x)}{\partial x}$ is the Jacobi matrix.
 Without loss of generality, we assume that each pair $(j,l)\in S$ is ordered with $j<l$.

Following the idea of~\cite{SL,Sussmann}, we introduce an extended system for~\eqref{Sigma}:
\begin{equation}
\dot x = \sum_{i=1}^m u_i f_i (x) + \sum_{(j,l)\in S} u_{jl}[f_j,f_l](x)\equiv \bar f(x,\bar u)
\label{Sigma_ext}
\end{equation}
with the control $\bar u = (u_1,u_2,...,u_m,{u_{jl}})'_{(j,l)\in S}$.
Because of the rank condition~\eqref{rank}, every smooth curve is a trajectory of system~\eqref{Sigma_ext}.
As subspaces spanned by the Lie brackets of vector fields $f_j(x)$ play a crucial role in the dynamics study of system~\eqref{Sigma}, we note that harmonic inputs naturally appear as optimal controls implementing the motion along a Lie bracket~\cite{Br1981,GJZ}.
A result on the convergence of solutions of system~\eqref{Sigma} to a solution of~\eqref{Sigma_ext} is established by H.J.~Sussmann and W.~Liu. It is shown in the paper~\cite{SL} that if a sequence of input functions $\{u^j(t)\}_{j=1}^\infty$ of class $L^1(0,\tau)$ satisfies certain boundedness condition and converges to an extented input $\bar u(t)$ in the iterated integrals sense, then solutions $x^j(t)$ of system~\eqref{Sigma} with initial data $x^j(0)=x^0$ converge to a solution $x^\infty(t)$ of system~\eqref{Sigma_ext}, uniformly with respect to $t\in [0,\tau]$. This result is stated for an extended system with higher order Lie brackets as well.
The problem of approximating a given trajectory of the extended system by trajectories of system~\eqref{Sigma} is solved in the paper~\cite{Liu} by using an unbounded sequence of oscillating controls with unbounded frequencies.
For a class of control systems with periodic solutions and small controls, an averaged control system is constructed in the paper~\cite{BP}. It is proved there that solutions of the averaged system approximate all solutions of the oscillating system as the frequency of oscillations tends to infinity.

In contrast to the above approach, we will use a time-varying feedback control $u=u(t,x)$ with bounded frequencies to implement certain decreasing condition for a Lyapunov functon along the trajectories of system~\eqref{Sigma}.
The rank condition~\eqref{rank} implies that any positive definite function $V(x)$ of class $C^1({\mathbb R}^n)$ may be taken as a control Lyapunov function for system~\eqref{Sigma_ext}, so its origin $x=0$ is stabilizable by a smooth feedback law $\bar u=\bar u(x)$, $\bar u(0)=0$.
Suppose that such a feedback $\bar u(x)$ is given, then our goal is to construct a time-varying feedback law $u=u(t,x)$ for the original system~\eqref{Sigma} in order to
approximate the flow of the closed-loop system~\eqref{Sigma_ext} in a suitable way. By exploiting this idea, we establish a result on the exponential stabilization in the sense of sampling controls and ``$\pi_\varepsilon$-solutions''.

We prove that, for systems satisfying the rank condition~\eqref{rank}, there exists a feedback $u=u^\varepsilon(t,x)$ such that any $\pi_\varepsilon$-solution $x(t)$ together with $u^\varepsilon(t,x(t))$ tend to zero exponentially, provided that $\varepsilon>0$ is small enough (Theorems~\ref{thm_sampling} and~\ref{thm_sampling2} in Section~\ref{sect_sampling}).
The proof of this result, given in Section~\ref{section_proof}, is based on Lyapunov's direct method and the representation of solutions by means of the Volterra series described in Section~\ref{section_oscillating}.
The construction of a stabilizing control $u=u^\varepsilon(t,x)$ is carried out explicitly in Section~\ref{section_Brockett} for the Brockett integrator. We show that such a feedback ensures exponential stability of the equilibrium.

\section{Stabilization with sampling controls}\label{sect_sampling}

For a given $\varepsilon>0$, we denote by $\pi_\varepsilon$ the {\em partition} of $[0,+\infty)$
into intervals
$$
I_j = [t_j,t_{j+1}),\;t_j=\varepsilon j,\quad j=0,1,2,\dots\; .
$$
The following definition extends the notion of ``$\pi$-trajectories'', introduced in~\cite{CLSS}, for the case of a time-varying feedback law.

{\bf Definition.} {\em
\label{definition_pi}
Assume given a feedback $u=h(t,x)$, $h:[0,+\infty)
\times D\to \mathbb R^m$, $\varepsilon>0$, and $x^0\in\mathbb R^n$.
A $\pi_\varepsilon$-solution of system~\eqref{Sigma} corresponding to $x^0\in D$ and  $h(t,x)$ is an absolutely continuous function $x(t)\in D$, defined for $t\in [0,+\infty)$, which satisfies the initial condition $x(0)=x^0$ and the following differential equations
$$
\dot x(t) = f(x(t),h(t,x(t_j))),\quad t\in I_j=[t_j,t_{j+1}),
$$
for each $j=0,1,2, \dots$ .
}

In order to stabilize system~\eqref{Sigma}, we will use a time-varying feedback control of the form
\begin{equation}
u^\varepsilon(t,x) = v(x) + \sum_{(i,l)\in S} a_{il}(x)\left\{\cos\left(\frac{2\pi k_{il}(x) }{\varepsilon}t\right)e_i+\sin\left(\frac{2\pi k_{il}(x) }{\varepsilon}t\right)e_l\right\}
\label{S1}
\end{equation}
on each interval $I_j$ of length $\varepsilon$, where $e_i$ denotes the $i$-th unit vector in $\mathbb R^m$, and functions $v(x)=\left(v_1(x),v_2(x),...,v_m(x)\right)'$, $a_{il}(x)$, $k_{il}(x)$ will be defined below.

Note that there is no control Lyapunov function {\em for the original system}~\eqref{Sigma} due to Artstein's theorem~\cite{Artstein} and Brockett's condition~\cite{Br1983}.
Even though a Lyapunov function may be constructed for system~\eqref{Sigma} in the sense of partial stability~\cite{Zuyev}, such partial formulation is not sufficient to establish an exponential stability result.
Because of the rank condition~\eqref{rank}, any differentiable positive definite function $V:D\to \mathbb R$ is a control Lyapunov function {\em for the extended system}~\eqref{Sigma_ext}.
Our main idea is to choose the feedback control~\eqref{S1} in order to approximate the direction of $-\nabla V(x)$ by trajectories of system~\eqref{Sigma}, where $\nabla V(x)$ is the gradient of $V(x)$. For this purpose, we fix $x\in D$ and $\varepsilon>0$, and consider the following system of second order algebraic equations
$$
\sum_{i=1}^m v_i f_i(x) + \frac{\varepsilon}{4\pi}\sum_{(i,j)\in S}\frac{a_{ij}^2}{k_{ij}}[f_i,f_j](x)+\frac{\varepsilon}{2}\sum_{i,j=1}^m v_i v_j\frac{\partial f_j(x)}{\partial x}f_i(x) +
$$
\begin{equation}
+\frac{\varepsilon}{2\pi}\sum_{i<j}\left(v_j\sum_{(q,i)\in S}\frac{a_{qi}}{k_{qi}}-v_i\sum_{(q,j)\in S}\frac{a_{qj}}{k_{qj}}\right)[f_i,f_j](x)= -\nabla V(x),
\label{Sigma1}
\end{equation}
with respect to variables $v_i$, $a_{ql}$, $i\in \{1,2,...,m\}$, $(q,l)\in S$,
 assuming that the numbers $k_{ql}\in \mathbb Z\setminus\{0\}$ are chosen without resonances, i.e.
\begin{equation}
|k_{ql}|\neq |k_{jr}|\quad \text{for\; all} \quad (q,l)\in S,(j,r)\in S, (q,l)\neq (j,r).
\label{nonres}
\end{equation}

Let us denote by $B_\rho(0)\subset \mathbb R^n$ the open ball of radius $\rho$ centered at $x=0$,
and let $\overline{B_\rho(0)}$ be its closure.
In this paper, we use the standard Euclidean norms for all vectors and treat $\frac{\partial^2 f_{ij}(x)}{\partial^2 x}$ as the Hessian matrix of the $j$-th component of $f_i(x)$.

The basic result of this paper is as follows.

\begin{theorem}
\label{thm_sampling}
Let $V(x)$ be a function of class $C^2(D)$ such that
\begin{equation}
 \|\nabla V(x)\|^2 \ge \alpha_1 V(x),\; V(x)\ge \beta_1 \|x\|^2,\;
V(0)=0,
\label{V_ineq1}
\end{equation}
and let
\begin{equation}
\left\|\frac{\partial f_i(x)}{\partial x}\right\|\le L, \quad \forall x\in D,\; i\in \{1,...,m\},
\label{Lip_L}
\end{equation}
with some positive constants $\alpha_1$, $\beta_1$, and $L$.
Assume that, for some $\rho_0>0$ and $\varepsilon_0>0$, algebraic system~\eqref{Sigma1} admits a solution
$$
 v_i=v^{\varepsilon}_i(x),\; a_{jl} = a^{\varepsilon}_{jl}(x),\; k_{jl}=k^{\varepsilon}_{jl}(x),\quad
 i\in \{1,...,m\},\,(j,l)\in S,
$$
defined for all $x\in \overline{B_{\rho_0}(0)}\subset D$ and $\varepsilon\in (0,\varepsilon_0]$,
such that condition~\eqref{nonres} holds and
\begin{equation}
\lim_{\varepsilon\to 0}\left(
\sup_{0<\|x\|\le \rho_0}\frac{\|v^{\varepsilon}(x)\|+\|a^{\varepsilon}(x)\|}{\|x\|^{1/3}}\varepsilon^{2/3}
\right)=0.
\label{a_v_asymptotics}
\end{equation}

Then there exist positive numbers $\rho\le \rho_0$ and $\bar \varepsilon\le \varepsilon_0$ such that,
for any $\varepsilon\in (0,\bar\varepsilon]$, there is a $\lambda=\lambda(\varepsilon)>0$:
\begin{equation}
x^0\in \overline{B_\rho(0)} \Rightarrow\; \|x(t)\|=O(e^{-\lambda t}),\; \|u^\varepsilon(t,x(t))\|
=O(e^{-\lambda t/3})\quad \text{as} \;\; t\to +\infty,
\label{exp_property}
\end{equation}
for each $\pi_\varepsilon$-solution $x(t)$ of system~\eqref{Sigma} with the control $u=u^\varepsilon(t,x)$ of form~\eqref{S1}.
\end{theorem}

Property~\eqref{exp_property} implies, in particular, that all $\pi_\varepsilon$-solutions $x(t)$ of the closed-loop system~\eqref{Sigma} and~\eqref{S1} with initial data $\|x^0\|\le \rho$ are defined for all $t\ge 0$.

To ensure the local solvability of equations~\eqref{Sigma1} in some $\Delta$-neigborhood of the point $x=0$,
we use the following lemma.

\begin{lemma}\label{lemma_solvability}
Assume that the vector fields $f_1(x)$, $f_2(x)$, ..., $f_m(x)$ satisfy
the rank condition~\eqref{rank} in a domain $D\subset \mathbb R^n$, $0\in D$, $|S|=n-m$,
and let $V\in C^2(D)$ be a positive definite function.
Then, for any small enough $\varepsilon>0$, there exists a $\Delta>0$ such that
algebraic system~\eqref{Sigma1} has a solution
$$
v^\varepsilon(x)=(v^\varepsilon_1(x),...,v^\varepsilon_m(x))',\; a^\varepsilon(x)=(a^\varepsilon_{jl}(x)_{(j,l)\in S})',\; k^\varepsilon(x)=(k^\varepsilon_{jl}(x)_{(j,l)\in S})',
$$
such that conditions~\eqref{nonres} hold for each $x\in B_\Delta(0)$. The above solution satisfies estimates
\begin{equation}
\|v^\varepsilon(x)\| \le M_v \|x\|,\; \|a^\varepsilon(x)\| \le M_a \sqrt{\frac{\|x\|}{\varepsilon}},\quad x\in B_\Delta(0),
\label{a_v_estimates}
\end{equation}
where positive constants $M_v$ and $M_a$ do not depend on $\varepsilon$.
\end{lemma}

The proof of Lemma~\ref{lemma_solvability} is based on the degree theory and will be presented in Section~\ref{section_proof}.
Lemma~\ref{lemma_solvability} allows us to formulate a local version of Theorem~\ref{thm_sampling} as follows.

\begin{theorem}
\label{thm_sampling2}
Assume that the vector fields $f_1(x)$, $f_2(x)$, ..., $f_m(x)$ satisfy the rank condition~\eqref{rank} with $|S|=n-m$ at $x=0$.
Then, for any positive definite quadratic form $V(x)$, there exist  constants
$\rho_0\ge \rho>0$ and $\varepsilon_0\ge \bar\varepsilon>0$ such that algebraic system~\eqref{Sigma1} admits a solution
$$
 v_i=v^{\varepsilon}_i(x),\; a_{jl} = a^{\varepsilon}_{jl}(x),\; k_{jl}=k^{\varepsilon}_{jl}(x),\quad  x\in \overline{B_{\rho_0}(0)}\subset D,\,\varepsilon\in (0,\varepsilon_0],
 $$
 $$
 i\in \{1,...,m\},\,(j,l)\in S,\,
$$
and, for any $\varepsilon\in (0,\bar\varepsilon]$, there is a $\lambda=\lambda(\varepsilon)>0$:
\begin{equation}
x^0\in \overline{B_\rho(0)} \Rightarrow\; \|x(t)\|=O(e^{-\lambda t}),\; \|u^\varepsilon(t,x(t))\|
=O(e^{-\lambda t/3})\quad \text{as} \;\; t\to +\infty,
\label{exp_property}
\end{equation}
for each $\pi_\varepsilon$-solution $x(t)$ of system~\eqref{Sigma} with the control $u=u^\varepsilon(t,x)$ of form~\eqref{S1}.
\end{theorem}
{\bf Proof.}
The assertion of Theorem~\ref{thm_sampling2} is a straightforward consequence of Theorem~\ref{thm_sampling}. To ensure condition~\eqref{a_v_asymptotics}, we use inequalities~\eqref{a_v_estimates} from Lemma~\ref{lemma_solvability}.
$\square$

The next section provides some technical results for the control design and stability analysis.
Then the proof of Theorem~\ref{thm_sampling} will be given in Section~\ref{section_proof}.

\section{Oscillating controls and representation of solutions}\label{section_oscillating}

Any solution $x(t)$ of system~\eqref{Sigma} with initial data $x(0)=x^0$ and controls $u_i=u_i(t)$, $u_i\in L^\infty[0,\tau]$ can be represented by means of the Volterra type series~(cf.~\cite{Boyd,NS}):
$$
x(\tau)=x^0+\sum_{i=1}^m f_i(x^0)\int_0^\tau u_i(t)dt + \frac{1}{2}\sum_{i,j=1}^m \frac{\partial f_j(x^0)}{\partial x}f_i(x^0)
\int_0^\tau u_i(t)dt \int_0^\tau u_j(t)dt+
$$
\begin{equation}
+\frac{1}{2}\sum_{i<j}[f_i,f_j](x^0)\int_0^\tau \int_0^t \left\{u_j(t)u_i(s)-u_i(t)u_j(s)\right\}ds\,dt + R(\tau).
\label{Volterra_series}
\end{equation}
Here, and in the sequel, $\frac{\partial f_j(x^0)}{\partial x}$ stands for the Jacobian matrix of $f_j(x)$ evaluated at $x=x^0$.
The remainder $R(\tau)$ of expansion~\eqref{Volterra_epsilon} is estimated by using the following lemma.

\begin{lemma}
\label{lemma_residual_Volterra}
Let $D\subset \mathbb R^n$ be a convex domain, and let $x(t)\in D$, $0\le t\le \tau$, be the solution of system~\eqref{Sigma} corresponding to an initial value $x(0)=x^0\in D$ and control $u\in C[0,\tau]$.
If the vector fields $f_1(x)$, $f_2(x)$, ..., $f_m(x)$ satisfy assumptions
 \begin{equation}
\left\|\frac{\partial f_i(x)}{\partial x}\right\|\le L,\; \left\|\frac{\partial^2 f_{ij}(x)}{\partial^2 x}\right\|\le H,\;\quad  i=\overline{1,m},\;j=\overline{1,n},
\label{Lip_constants}
\end{equation}
 in $D$ with some constants $H$, $L>0$,
then the remainder $R(\tau)$ of the Volterra expansion~\eqref{Volterra_series} satisfies the following estimate:
$$
\|R(\tau)\|\le \frac{M}{L}\left\{e^{LU\tau}-\frac{1}{2}\left((LU\tau+1)^2+1\right)\right\}+
$$
\begin{equation}
+\frac{HM^2\sqrt{n}}{4L^3}\left\{\left(e^{LU\tau}-2\right)^2+2LU\tau -1\right\}=\frac{M( L^2+HM\sqrt{n})}{6}U^3\tau^3+O(U^4\tau^4).
\label{remainder}
\end{equation}
Here
$$
M = \max_{1\le i\le m} \|f_i(x^0)\|,\; U = \max_{0\le t\le \tau} \sum_{i=1}^m |u_i(t) |.
$$
\end{lemma}
The proof of Lemma~\ref{lemma_residual_Volterra} is given in Section~\ref{section_proof}.

In order to use the control strategy~\eqref{S1}, we consider a family of open-loop controls
\begin{equation}
u_i(t) = v_i + \sum_{(j,l)\in S} a_{jl}\left\{\delta_{ij}\cos\left(\frac{2\pi k_{jl} }{\varepsilon}t\right)+\delta_{il}\sin\left(\frac{2\pi k_{jl} }{\varepsilon}t\right)\right\},\;i=1,2,...,m,
\label{controls_m}
\end{equation}
depending on parameters $v=(v_1,v_2,...,v_m)'\in\mathbb R^m$, $a=\left(a_{jl}\right)'_{(j,l)\in S}\in \mathbb R^{n-m}$, $k=\left(k_{jl}\right)'_{(j,l)\in S}\in (\mathbb Z\setminus \{0\})^{n-m}$, and $\varepsilon>0$. Here $\delta_{ij}$ is the Kronecker delta.

By computing the integrals in~\eqref{Volterra_series} for functions $u_i=u_i(t)$ given by~\eqref{controls_m} and exploiting assumption~\eqref{nonres}, we get
$$
x(\varepsilon)=x^0 + \varepsilon \sum_{i=1}^m v_i f_i(x^0)
+ \frac{\varepsilon^2}{2}\sum_{i,j=1}^m v_i v_j\frac{\partial f_j(x^0)}{\partial x}f_i(x^0)+
$$
\begin{equation}
+\frac{\varepsilon^2}{4\pi}\sum_{i<j}[f_i,f_j](x^0)\sum_{(q,l)\in S}\frac{a_{ql}}{k_{ql}}\left\{\delta_{jl}(a_{ql}\delta_{iq}-2v_i)-\delta_{il}(a_{ql}\delta_{jq}-2v_j)\right\}+R(\varepsilon).
\label{Volterra_epsilon}
\end{equation}

To estimate the decay rate of the function $V(x(t))$, we use the following lemma.

\begin{lemma}
\label{Lyapunov_decay}
Let $V(x)$ be a function of class $C^2(D)$ such that inequalities
\begin{equation}
\beta_1 \|x\|^2 \le V(x) \le \beta_2 \|x\|^2,\quad \beta_1>0,
\label{V_ineq}
\end{equation}
\begin{equation}
\alpha_1 V(x) \le \|\nabla V(x)\|^2\le \alpha_2 V(x),\quad \alpha_1>0,
\label{nabla_ineq}
\end{equation}
\begin{equation}
\left\|\frac{\partial^2 V(x)}{\partial x^2}\right\|\le \mu
\label{Hesse_ineq}
\end{equation}
hold with some constants $\alpha_1$, $\alpha_2$, $\beta_1$, $\beta_2$, $\mu$ in a convex domain $D\subset \mathbb R^n$.

If $x:[0,\varepsilon]\to D$ is a function such that
\begin{equation}
x(\varepsilon) = x(0) - \varepsilon \nabla V(x(0)) +r_\varepsilon,\; x(0)\neq 0,
\label{x_variation}
\end{equation}
with some $r_\varepsilon\in \mathbb R^n$, then
\begin{equation}
V(x(\varepsilon))\le V(x(0))\left\{1-\alpha_1\varepsilon + \frac{\alpha_2 \varepsilon^2 \mu}{2}+\frac{\mu\|r_\varepsilon\|^2}{2\beta_1 \|x(0)\|^2}+\frac{\sqrt{\alpha_2}(1+\varepsilon\mu)\|r_\varepsilon\|}{\sqrt{\beta_1}\|x(0)\|}\right\}.
\label{V_estimate}
\end{equation}
\end{lemma}
{\bf Proof.}
Let us denote $x^0=x(0)$, $y=-\varepsilon \nabla V(x^0)+r_\varepsilon$, and apply Taylor's theorem for the function $V(x^0+y)$ with the Lagrange form of the remainder:
\begin{equation}
V(x^0+y)=V(x^0)+\sum_{i=1}^n\left.\frac{\partial V(x)}{\partial x_i}y_i\right|_{x=x^0} +\frac{1}{2}\sum_{i,j=1}^n\left.\frac{\partial^2 V(x)}{\partial x_i\partial x_j}y_iy_j\right|_{x=x^0+\theta y},
\label{taylor_v}
\end{equation}
where $\theta\in (0,1)$. By applying the Cauchy--Schwartz inequality to expansion~\eqref{taylor_v} and exploiting assumptions~\eqref{nabla_ineq}, \eqref{Hesse_ineq}, we get the following estimate:
$$
V(x(\varepsilon))\le V(x^0)-\varepsilon\left(1 - \frac{\varepsilon \mu }{2}\right)\|\nabla V(x^0)\|^2 +(1+\varepsilon\mu)\|\nabla V(x^0)\|\,\|r_\varepsilon\|+\frac{\mu}{2}\|r_\varepsilon\|^2\le
$$
$$
\le \left(1-\alpha_1\varepsilon + \frac{\alpha_2 \varepsilon^2 \mu}{2}\right)V(x^0) + (1+\varepsilon\mu)\sqrt{\alpha_2 V(x^0)}\|r_\varepsilon\| + \frac{\mu \|r_\varepsilon\|^2}{2}\le
$$
\begin{equation}
\le V(x^0)\left(1-\alpha_1\varepsilon + \frac{\alpha_2 \varepsilon^2 \mu}{2}+\frac{\sqrt{\alpha_2}(1+\varepsilon\mu)\|r_\varepsilon\|}{\sqrt{V(x^0)}}+\frac{\mu \|r_\varepsilon\|^2}{2V(x^0)}\right)
\label{decreasing_est}
\end{equation}
if $V(x^0)\neq 0$. Then the application of estimate~\eqref{V_ineq} to~\eqref{decreasing_est} yields ineqiality~\eqref{V_estimate}.
$\square$

By using Lemmas~\ref{lemma_residual_Volterra} and~\ref{Lyapunov_decay} for $\pi_\varepsilon$-solutions of system~\eqref{Sigma} corresponding to a partition $\pi_\varepsilon=\{\varepsilon j\}_{j\ge 0}$ and control $u=u^\varepsilon(t,x)$, we prove Theorem~\ref{thm_sampling}.

\section{Proof of the main result}\label{section_proof}

In order to prove Theorem~\ref{thm_sampling}, let us first prove auxiliary lemmas.

\textbf{Proof of Lemma~\ref{lemma_solvability}.}
Let us enumerate the elements of $S$ in~\eqref{rank} as
$$
S=\{(i_1,j_1), \; (i_2,j_2),\; ...\, , (i_{n-m},j_{n-m})\}
$$
and introduce the $n\times n$-matrix
$$
A(x)=\bigl(f_1(x),...,f_m(x),[f_{i_1},f_{j_1}](x),...,[f_{i_{n-m}},f_{j_{n-m}}(x)]\bigr).
$$
As the vector fields $f_i(x)$ satisfy the rank condition~\eqref{rank}, there is a closed bounded domain $\Omega\subset D$, $0\in {\rm int}\, D$ such that the map
\begin{equation}
\Phi(x)=-A^{-1}(x)\nabla V(x)
\label{phi_def}
\end{equation}
is defined for each $x\in \Omega$ and continuous.
To study the solvability of equations~\eqref{Sigma1}, we introduce new variables
$$
\tilde a_{ij} = \frac{a^2_{ij}}{4\pi k_{ij}},\quad (i,j)\in S,
$$
and rewrite system~\eqref{Sigma1} as follows
\begin{equation}
F_x(\xi)+G_x(\xi) = 0,\quad \xi=\left(v_1,...,v_m,\varepsilon\tilde a_{i_1 j_1},...,\varepsilon\tilde a_{i_{n-m} j_{n-m}}\right)',
\label{alg_sys_symbolic}
\end{equation}
where
$$
F_x(\xi) = \xi -\Phi(x),
$$
$$
\frac{2}{\varepsilon}G_x(\xi) =  \sum_{i,j=1}^m v_i v_jA^{-1}(x)\frac{\partial f_j(x)}{\partial x}f_i(x)
+
$$
$$
+\frac{2}{\sqrt{\pi}}\sum_{i<j}\left(v_j\sum_{(q,i)\in S}\sqrt{\frac{|\tilde a_{qi}|}{\bar k_{qi}}}-v_i\sum_{(q,j)\in S}\sqrt{\frac{|\tilde a_{qj}|}{\bar k_{qj}}}\right)A^{-1}(x)[f_i,f_j](x).
$$
Here the integer constants $\bar k_{ij}$ may be chosen as
$$
\bar k_{i_1 j_1}=1,\; \bar k_{i_2 j_2}=2,\; ...,\; \bar k_{i_{n-m} j_{n-m}}=n-m.
$$
If $\xi=\left(v_1,...,v_m,\varepsilon\tilde a_{i_1 j_1},...,\varepsilon\tilde a_{i_{n-m} j_{n-m}}\right)'$ is a solution of system~\eqref{alg_sys_symbolic} for given $x\in\mathbb R^n$ and $\varepsilon>0$, then the components of a solution of equations~\eqref{Sigma1} are
\begin{equation}
v_1, v_2, ..., v_m, \; a_{ij} = 2\sqrt{\pi\bar k_{ij}|\tilde a_{ij}|}\,{\rm sign }\,\tilde a_{ij},
\; (i,j)\in S,
\label{solutions}
\end{equation}
with
$k_{ij} = \bar k_{ij}$ if $\tilde a_{ij}\ge 0$ and $k_{ij}=-\bar k_{ij}$ otherwise.
Thus, the solvability of system~\eqref{Sigma1} is reduced to the study of equation~\eqref{alg_sys_symbolic}.

Note that $\Phi(x)=0$ for $x=0$ as $V(x)$ is positive definite,
so $\xi=0$ is a solution of equation~\eqref{alg_sys_symbolic} for $x=0$.

To prove the existence of solutions for equation~\eqref{alg_sys_symbolic},
we find a $\Delta>0$ and show that the degree of a continuous map
\begin{equation}
\xi\in S_\rho \mapsto \frac{F_x(\xi)+G_x(\xi)}{\|F_x(\xi)+G_x(\xi)\|}\in S_1,
\label{map_sphere}
\end{equation}
is equal to 1, under a suitable choice of $\rho>0$ depending on $x$ if $0< \|x\|<\Delta$, $\overline{B_\Delta(0)}\subset \Omega$.
Here the spheres
$$S_\rho = \{\xi\in\mathbb R^n\,|\, \|\xi\|=\rho\}\;\text{and}\;
S_1 = \{\xi\in\mathbb R^n\,|\, \|\xi\|=1\}$$
are oriented as $(n-1)$-spheres in $\mathbb R^n$.

As $\Omega$ is compact then there exist positive constants $M_0$, $M_1$, and $L$ such that
$$
\|A^{-1}(x)\|\le M_0,\;\|f_i(x)\|\le M_1,\; \left\|\frac{\partial f_i(x)}{\partial x}\right\|\le L,\;i=1,2,...,m,\quad \forall x\in\Omega.
$$
If $\|\xi\|\in S_\rho$ then the Cauchy--Schwartz and triangle inequalities yield
$$
\|G_x(\xi)\|\le {\varepsilon \bar M}\sum_{i\le j}|v_i v_j|+\frac{2\varepsilon \bar M}{\sqrt{\pi}}(n-m)^{3/4}\|\tilde a\|^{1/2}\sum_{i<j}\left(|v_i|+|v_j|\right)<
$$
\begin{equation}
< {\varepsilon \bar M n}\|v\|^2 + \frac{4\varepsilon \bar M n^{9/4}}{\sqrt{\pi}}\|v\|\,\|\tilde a\|^{1/2}\le \bar M n \sqrt{\varepsilon\rho^3}\left\{\frac{4n^{5/4}}{\sqrt{\pi}}+\sqrt{\varepsilon \rho}\right\},
\label{G_estimate}
\end{equation}
where $\bar M=L M_0 M_1$.
We have also exploited formula~\eqref{solutions} together with the following properties of components of $\xi$ here:
\begin{equation}
\|v\|\le \|\xi\|,\;\|\tilde a\| \le \frac{1}{\varepsilon}\|\xi\|.
\label{a_v_xi}
\end{equation}

Then the maps $F_x(\xi)+G_x(\xi)$ and $\xi$ are homotopic on $S_\rho$ provided that
\begin{equation}
\|\xi\|> \|\Phi(x)\|+\|G_x(\xi)\|,\quad \forall \xi\in S_\rho.
\label{homotopy_condition}
\end{equation}
To satisfy condition~\eqref{homotopy_condition}, we observe that $\|\nabla V(x)\|=O(\|x\|)$ in a neighborhood of $x=0$ for a positive definite function $V(x)$.
Hence, there exist positive constants $\bar \Delta$ and $\psi$ such that
\begin{equation}
\|\Phi(x)\| \le \psi \|x\|,\quad \forall x\in B_{\bar\Delta}(0).
\label{phi_estimate}
\end{equation}
By taking into account inequalities~\eqref{G_estimate} and~\eqref{phi_estimate}, we conclude that condition~\eqref{homotopy_condition} is satisfied for $x\in B_{\bar\Delta}(0)$ if
$$
\rho\ge \psi\|x\| + \bar M n \sqrt{\varepsilon\rho^3}\left\{\frac{4n^{5/4}}{\sqrt{\pi}}+\sqrt{\varepsilon \rho}\right\}.
$$
The function $\phi(\rho)=\rho-\bar M n \sqrt{\varepsilon\rho^3}\left\{\frac{4n^{5/4}}{\sqrt{\pi}}+\sqrt{\varepsilon \rho}\right\}=\rho+O(\rho^{3/2})$ takes positive values for $\rho\in(0,\rho_0)$, where
$$
\sqrt{\rho_0}=r_0 = \frac{2\pi n^{5/4}}{\sqrt{\varepsilon\pi}}\left(\sqrt{1+\frac{\pi}{4\bar Mn^{7/2}}}-1\right)
$$ is the positive root of the following equation
$$
r_0^2 + \frac{4n^{5/4}}{\sqrt{\varepsilon\pi}}r_0-\frac{1}{\bar Mn\varepsilon}=0.
$$
Then we choose any $\bar r\in (0,r_0)$ and check that
$$
\phi(\rho) \ge K \rho \quad \text{for}\; 0\le \rho\le {\bar r}^2,
$$
where
$$
K = 1-\sqrt{\varepsilon}\bar Mn{\bar r}^3\left\{\frac{4n^{5/4}}{\sqrt{\pi}}+\bar r\sqrt{\varepsilon}\right\}>0.
$$
Let us take
$$
\Delta ={\rm min}\left\{\bar\Delta,\frac{K{\bar r}^2}{\psi}\right\}>0
$$
and observe that, for any $x:0< \|x\|<\Delta$, condition~\eqref{homotopy_condition} holds if
\begin{equation}
\rho=\frac{\psi}{K}\|x\|.
\label{rho_formula}
\end{equation}
Homotopic equivalence of the maps $F_x(\xi)+G_x(\xi)$ and $\xi$ on $S_\rho$, ensured by~\eqref{homotopy_condition}, implies that the degree of the map~\eqref{map_sphere} is equal to 1, i.e. to the degree of the map $\frac{\xi}{\|\xi\|}:S_\rho\to S_1$.
By exploiting the degree principle (see, e.g.,~\cite{KZ}), we conclude that there exists a $\xi\in B_\rho(0)$ such that $F_x(\xi)+G_x(\xi)=0$, which means the existence of a solution to equation~\eqref{Sigma1} according to formulas~\eqref{solutions} if $x\in B_\Delta(0)$.
Then estimates~\eqref{a_v_estimates} follow from inequalities~\eqref{a_v_xi} and~\eqref{rho_formula}.
$\square$

A useful a priori estimate of the solutions of system~\eqref{Sigma} is given by the following lemma.

\begin{lemma}\label{lemma_apriori}
Let $x(t)\in D\subset \mathbb R^n$, $0\le t\le \tau$, be a solution of system~\eqref{Sigma} with the control $u\in C[0,\tau]$, and let
$$
\|f_i(x')-f_i(x'')\|\le L \|x'-x''\|,\quad \forall x',x''\in D,\; i= 1,2,..., m.
$$
Then
\begin{equation}
\|x(t)-x(0)\|\le \frac{M}{L}(e^{LUt}-1),\quad t\in [0,\tau],
\label{apriori_est}
\end{equation}
where
$$
M=\max_{1\le i\le m}\|f_i(x(0))\|,\; U= \max_{0\le t\le \tau}\sum_{i=1}^m |u_i(t)|.
$$
\end{lemma}
{\bf Proof.}
By differentiating the function $w(t)=\|x(t)-x(0)\|$ along the trajectory of system~\eqref{Sigma},
we get
$$
\frac{d}{dt}w^2(t)=2\left(x(t)-x(0),\sum_{i=1}^m u_i(t)f_i(x(t))\right)\le
$$
$$
\le 2U w(t)\max_{1\le i\le m}\|f_i(x(t))-f_i(x(0))+f_i(x(0))\|\le 2U w(t)(Lw(t)+M),
$$
so,
\begin{equation}
\dot w(t)\le U (Lw(t)+M),\quad t>0.
\label{wdot}
\end{equation}
We solve the comparison equation for differential inequality~\eqref{wdot} to obtain the following estimate~(cf.~\cite[Chap.~III]{Hartman}):
$$
w(t)\le \frac{M}{L}(e^{LUt}-1),\quad t\in [0,\tau].
$$
This proves estimate~\eqref{apriori_est}.
$\square$

Now we use Lemma~\ref{lemma_apriori} to prove Lemma~\ref{lemma_residual_Volterra}.

{\em Proof of Lemma~\ref{lemma_residual_Volterra}.}
For a solution $x(t)$ of differential equation~\eqref{Sigma} with the initial condition $x(0)=x^0$ and control $u\in C[0,\tau]$, we represent the coordinates of $\Delta x(t)=x(t)-x^0$ by the following integral equations:
$$
\Delta x_k(\tau)=\sum_{i=1}^m \int_0^\tau u_i(t)f_{ik}\left(x^0+\sum_{j=1}^m \int_0^t u_j(s)f_j(x^0+\Delta x(s))ds\right)dt=
$$
$$
=\sum_{i=1}^m \int_0^\tau u_i(t)\left\{ f_{ik}(x^0)+\frac{\partial f_{ik}(x^0)}{\partial x}\sum_{j=1}^m \int_0^t u_j(s)\left(f_j(x^0)+\frac{\partial f_j(\xi(s))}{\partial x}\Delta x(s)\right)ds\right.+
$$
\begin{equation}
+\left. \frac{1}{2}\left(\frac{\partial^2 f_{ik}(\eta(t))}{\partial x^2}\Delta x(t),\Delta x(t)\right)\right\}dt,\quad k=1,2,..., n.
\label{integral_eq}
\end{equation}
Expression~\eqref{integral_eq} is obtained by Taylor's theorem with the Lagrange form of the remainder,
$0\le \|\xi(s)-x^0\|\le \|\Delta x(s)\|$, $0\le \|\eta(t)-x^0\|\le \|\Delta x(t)\|$, for $0\le s\le t\le \tau$. Comparing formula~\eqref{Volterra_series} with~\eqref{integral_eq}, we get
$$
R_k(\tau) = \sum_{i,j=1}^m \frac{\partial f_{ik}(x^0)}{\partial x}\int_0^\tau \int_0^t u_i(t)u_j(s)
\frac{\partial f_j(\xi(s))}{\partial x}\Delta x(s)\,ds\, dt+
$$
\begin{equation}
+\frac{1}{2}\sum_{i=1}^m \int_0^\tau \left(\frac{\partial^2 f_{ik}(\eta(t))}{\partial x^2}\Delta x(t),\Delta x(t)\right)u_i(t)\, dt.
\label{remainder_int}
\end{equation}
We use estimate~\eqref{apriori_est} from Lemma~\ref{lemma_apriori} and the triangle inequality together with the Cauchy--Schwartz inequality to evaluate the Euclidean norm of the vector $R(\tau)=(R_1(\tau),...,R_n(\tau))'$ in~\eqref{remainder_int}:
$$
\|R(\tau)\| \le L^2 U^2 \int_0^\tau \int_0^t \|\Delta x(s)\|\,ds\, dt + \frac{H U \sqrt{n}}{2}\int_0^\tau \|\Delta x(t)\|^2 dt \le
$$
\begin{equation}
\le \frac{M}{L}\left\{e^{LU\tau} -\frac{1}{2}\left((LU\tau+1)^2+1\right)\right\}+
\frac{H M^2 \sqrt{n}}{4 L^3}\left\{ \left(e^{LU\tau}-2\right)^2 +2LU\tau -1 \right\}.
\label{remainder_norm}
\end{equation}
The right-hand side of formula~\eqref{remainder} is obtained as the Taylor expansion of formula~\eqref{remainder_norm} with respect to $U\tau$. $\square$

{\em Proof of Theorem~\ref{thm_sampling}}.
Let us denote $D_0=\overline{B_{\rho_0}(0)}\subset D$ and choose a positive number $\hat\varepsilon\le \varepsilon_0$ such that all solutions $x(t)$ of system~\eqref{Sigma} are well defined on $t\in [0,\varepsilon]$ for each $\varepsilon\in (0,\hat\varepsilon]$, provided that $x^0=x(0)\in D_0$ and the control $u=u^\varepsilon(t,x^0)$ is given by formula~\eqref{S1} with parameters
$$v_i=v_i^\varepsilon(x^0),\; a_{jl}=a_{jl}^\varepsilon(x^0),\;k_{jl}=k_{jl}^\varepsilon(x^0)
$$
obtained from algebraic system~\eqref{Sigma1}.

We define
\begin{equation}
M = \sup_{x\in D_0,\, 1\le i \le m} \|f_i(x)\|,\; d=\inf_{x\in D_0,\, y\in \partial D} \|x-y\|>0.
\label{M_const}
\end{equation}
If $D={\mathbb R}^n$ then we take $d=+\infty$ and $\hat\varepsilon=\varepsilon_0$, otherwise $\hat\varepsilon\le \varepsilon_0$ is obtained as a positive solution of the inequality
\begin{equation}
\frac{M}{L}\left(e^{L U(\hat\varepsilon)\hat\varepsilon}-1\right)<d,
\label{ineq_d}
\end{equation}
where $L$ is given in condition~\eqref{Lip_L},
$$
U(\hat\varepsilon) = \sup_{t\in [0,\hat\varepsilon],\, x\in D_0}\sum_{i=1}^m |u_i^{\hat\varepsilon}(t,x)|,
$$
and $u_i^{\hat\varepsilon}(t,x)$ are given by~\eqref{S1}.
Condition~\eqref{a_v_asymptotics} implies that $U(\hat\varepsilon)\hat\varepsilon\to 0$ as $\hat\varepsilon\to 0$. Thus, the set of solutions $\hat\varepsilon\in (0,\varepsilon_0]$ of inequality~\eqref{ineq_d} is not empty. Let $\hat\varepsilon$ be such a solution, then from inequality~\eqref{ineq_d} and Lemma~\ref{lemma_apriori} it follows that
\begin{equation}
\|x(t)-x^0\|<d,\quad t\in [0,\varepsilon],
\label{ineq_d2}
\end{equation}
for each solution $x(t)$ of system~\eqref{Sigma} with $x^0\in D_0$ and $u=u^\varepsilon(t,x^0)$,
$\varepsilon\in (0,\hat\varepsilon]$. Inequality~\eqref{ineq_d2} means that $x(t)\in D$ for $t\in [0,\varepsilon]$.

Let $V(x)$ be a function that satisfies conditions~\eqref{V_ineq1}.
We introduce level sets
$$
L_c = \{x\in D\,|\, V(x)\le c\}
$$
and define
 $$
c_0 = \inf_{x\in D\setminus D_0}V(x),\; \rho = \inf_{x\in D\setminus L_{c_0}}\|x\|.
 $$
 It is easy to see that $c_0$ and $\rho\le \rho_0$ are positive numbers as $V(x)$ is positive definite.
 By the construction,
 $$
 \overline{B_\rho(0)}\subseteq L_{c_0}\subseteq D_0\;\;\text{and}\;\;
 L_c \subseteq L_{c_0}
 $$
 for each $c\le c_0$.

The next step is to show that, if $\varepsilon>0$ is small enough, then there exists a positive $\sigma=\sigma(\varepsilon)<1$ such that
\begin{equation}
V(x(\varepsilon))\le \left(1-\sigma\right)V(x^0),
\label{v_decreasing}
\end{equation}
for any solution $x(t)$ of system~\eqref{Sigma} with the initial data $x^0\in L_{c_0}$ and the control $u=u^\varepsilon(t,x^0)$ given by~\eqref{S1}.

As $V\in C^2(D)$ is positive definite then $\nabla V(0)=0$, and Taylor's theorem implies the following inequality:
\begin{equation}
V(x)\le \beta_2 \|x\|^2,\quad \forall x\in D_0,
\label{V_ineq2}
\end{equation}
where
$$
2 \beta_2 = \mu = \sup_{x\in D_0} \left\|\frac{\partial^2 V(x)}{\partial x^2}\right\|
$$
is finite by Weierstrass's theorem due to the compactness of $D_0$.
By applying similar argumentation to the function $\|\nabla V(x)\|^2$, we conclude that
$$
\|\nabla V(x)\|^2 \le \bar \alpha_2\|x\|^2,\quad \forall x\in D_0,
$$
with some positive constant $\bar \alpha_2$. Because of conditions~\eqref{V_ineq1}, it follows that
\begin{equation}
\|\nabla V(x)\|^2 \le \alpha_2 V(x),\quad \forall x\in D_0,
\label{V_ineq3}
\end{equation}
where $\alpha_2 = \bar \alpha_2 / \beta_1 >0$.
Inequalities~\eqref{V_ineq1}, \eqref{V_ineq2}, and \eqref{V_ineq3} imply that all conditions of
Lemma~\ref{Lyapunov_decay} are satisfied in $D_0$ if $x(t)$ ($0\le t\le \varepsilon$) is a solution of system~\eqref{Sigma} with the control $u=u^\varepsilon(t,x^0)$, $x^0\in D_0$.

In order to satisfy condition~\eqref{v_decreasing}, it suffices to assume that
\begin{equation}
\alpha_1 \varepsilon - \frac{\alpha_2\varepsilon^2\mu}{2}-\frac{\mu \|R(\varepsilon)\|^2}{2\beta_1\|x^0\|^2}-\frac{\sqrt{\alpha_2}(1+\varepsilon\mu)\|R(\varepsilon)\|}{\sqrt{\beta_1}\|x^0\|}\ge \sigma,\;\forall x^0\in D_0\setminus\{0\}
\label{sufficient_decreasing}
\end{equation}
because of Lemma~\ref{Lyapunov_decay}.
Here the remainder $R(\varepsilon)$ of the Volterra series can be estimated by
Lemma~\ref{lemma_residual_Volterra} as follows:
\begin{equation}
\|R(\varepsilon)\| \le \bar H W^3(x^0)\varepsilon^3\quad \text{for}\;\; x^0\in D_0\; \text{and}\; W(x^0)\varepsilon\le 1.
\label{R_est}
\end{equation}
Here $\bar H$ is a positive constant,
\begin{equation}
W(x^0) = \sup_{t\in [0,\varepsilon]}\sum_{i=1}^m |u_i^{\varepsilon}(t,x^0)|,
\label{W_def}
\end{equation}
and $u^{\varepsilon}(t,x^0)$ is given by~\eqref{S1}. Condition~\eqref{a_v_asymptotics} together with representation~\eqref{S1} implies that
\begin{equation}
W(x^0) \le C \|x^0\|^{1/3} \varepsilon^{-2/3}
\label{W_est}
\end{equation}
with some positive constant $C$ for all $x^0\in D_0$.

Estimates~\eqref{R_est} and~\eqref{W_est} imply that condition~\eqref{sufficient_decreasing} holds if
\begin{equation}
\alpha_1  - \frac{\alpha_2\mu\varepsilon}{2}-\frac{\mu \bar H^2 C^6\varepsilon^2}{2\beta_1}-\frac{\sqrt{\alpha_2}(1+\varepsilon\mu)\bar H C^3 \varepsilon}{\sqrt{\beta_1}}\ge \bar\sigma.
\label{decreasing2}
\end{equation}
Here $\bar \sigma=\sigma/\varepsilon$ is a positive number. As $\alpha_1$ is positive, we conclude that there exist $\varepsilon_{max}>0$ and $\bar\sigma>0$ such that inequality~\eqref{decreasing2} holds for all $\varepsilon\in (0,\varepsilon_{max}]$.
Without loss of generality we suppose that such $\varepsilon_{max}$ corresponds to the assumption of
formula~\eqref{R_est}, i.e. $W(x^0)\varepsilon_{max}\le 1$ for all $x^0\in D_0$.
Thus we have proved that condition~\eqref{sufficient_decreasing} is satisfied for each $\varepsilon\in (0,\varepsilon_{max}]$ with $\sigma =\sigma(\varepsilon) = \min(\bar \sigma\varepsilon,1)$.
Let us define $\bar\varepsilon=\min(\hat\varepsilon,\varepsilon_{max})$, where $\hat\varepsilon$ is a positive solution of inequality~\eqref{ineq_d}. Then inequality~\eqref{v_decreasing} holds for any
$\varepsilon\in (0,\bar\varepsilon]$ with $\sigma=\sigma(\varepsilon)\le 1$ provided that $x^0\in L_{c^0}$.

If $x^0\in B_\rho(0)$, $\varepsilon\in (0,\bar\varepsilon]$ and $u^\varepsilon(t,x)$ is given by formula~\eqref{S1}, then the corresponding $\pi_\varepsilon$-solution of system~\eqref{Sigma} $x(t)$ is well-defined: $$x(n\varepsilon)\in L_{c_0}\subseteq D
_0\quad \text{ for}\;\; n=0,1,2, ... ,
 $$
 and $x(t)\in D$ for all  $t\ge 0$ because of inequality~\eqref{ineq_d2}.
By iterating inequality~\eqref{v_decreasing} for $x^0\in \overline{B_\rho(0)}\subseteq L_{c_0}$, we conclude that
\begin{equation}
\|x(t)\|\le \sqrt{\frac{\beta_2}{\beta_1}}\|x^0\|e^{-\bar\lambda t} \quad \text{for}\;\; t=0,\varepsilon, 2\varepsilon, ... ,
\label{v_est_sampling}
\end{equation}
where
$$
\bar\lambda = -\frac{{\rm ln}(1-\sigma)}{2\varepsilon}>0\quad \text{if }\;\;\sigma<1,
$$
and $\bar\lambda$ is an arbitrary positive number if $\sigma=1$.
For an arbitrary $t\ge 0$, we denote the integer part of $\frac{t}{\varepsilon}$ as $N=\left[\frac{t}{\varepsilon}\right]$
and denote $\tau = t-N\varepsilon<\varepsilon$.
Then we apply inequality~\eqref{v_est_sampling} together with Lemma~\ref{lemma_apriori} to estimate $x(t)$:
$$
\|x(t)\|=\|x(t)-x(N\varepsilon)+x(N\varepsilon)\|\le \|x(N\varepsilon)\|+\|x(t)-x(N\varepsilon)\|\le
$$
\begin{equation}
\le \sqrt{\frac{\beta_2}{\beta_1}}\|x^0\|e^{-\bar\lambda N\varepsilon}+
\frac{M}{L}\left(e^{L W(x(N\varepsilon))\varepsilon}-1\right),
\label{x_norm}
\end{equation}
where  $L$, $M$, and $W(x)$ are defined in~\eqref{Lip_L}, \eqref{M_const}, and~\eqref{W_def}, respectively. Estimates~\eqref{W_est} and~\eqref{v_est_sampling} imply the following asymptotic representation:
$$
W(x(N\varepsilon))=O(\|x(N\varepsilon)\|^{1/3})=O(e^{-\bar\lambda N\varepsilon/3})\quad\text{as}\;\;N\to +\infty.
$$
Then if follows from inequality~\eqref{x_norm} that
$$
\|x(t)\| = O(e^{-\lambda t})\quad \text{as}\;\; t\to +\infty
$$
with $\lambda=\bar\lambda/3>0$. By using formulas~\eqref{W_def} and~\eqref{W_est}, we conclude that
$$
\|u^\varepsilon(t,x(t))\|=O(e^{-\lambda t/3})\quad \text{as}\;\; t\to +\infty.
$$
$\square$

\section{Stabilization of the Brockett integrator}\label{section_Brockett}

Consider the control system known as the Brockett integrator~\cite{Br1983}:
\begin{equation}
\dot x_1= u_1,\; \dot x_2=u_2,\; \dot x_3 = u_1 x_2 - u_2 x_1,
\label{Brockett3}
\end{equation}
where $x=(x_1,x_2,x_3)'\in {\mathbb R}^3$ is the state and $u=(u_1,u_2)'\in {\mathbb R}^2$ is the control.
The stabilization problem for system~\eqref{Brockett3} has been has been the subject of many publications over the past three decades (see, e.g., the book~\cite{Bloch} and references therein). In particular, it is shown that system~\eqref{Brockett3} can be exponentially stabilized by a time-invariant feedback law for the initial values in some open and dense set $\Omega\neq \mathbb R^3$, $0\notin {\rm int}\,\Omega$~\cite{Astolfi}. In this section, we construct a time-varying feedback law explicitly in order to stabilize system~\eqref{Brockett3} exponentially for all initial data.

System~\eqref{Brockett3} satisfies the rank condition of form~\eqref{rank} with $S=\{(1,2)\}$,
$$
{\rm span} \{f_1(x),f_2(x),[f_1,f_2](x)\}={\mathbb R}^3\quad \text{for each}\;\; x\in\mathbb R^3,
$$
where the vector fields are $f_1(x)=(1,0,x_2)'$, $f_2(x)=(0,1,-x_1)'$,
$$
[f_1,f_2](x) = \frac{\partial{f_2(x)}}{\partial{x}}f_1(x) - \frac{\partial{f_1(x)}}{\partial{x}}f_2(x) =(0,0,-2)'.
$$
The family of controls~\eqref{controls_m} takes the form
\begin{equation}
\begin{array}{l}
u_1 (t) = v_1 + a_{12}\cos\left(\frac{2\pi k_{12}}{\varepsilon}t \right),\\
u_2 (t) = v_2 + a_{12}\sin\left(\frac{2\pi k_{12}}{\varepsilon}t \right),\quad k_{12}\in\mathbb Z\setminus \{0\}.
\label{controls}
\end{array}
\end{equation}

For an arbitrary initial value $x^0=(x_1^0,x_2^0,x_3^0)'\in {\mathbb R}^3$ at $t=0$, the solution $x(t)$ of system~\eqref{Brockett3} with controls~\eqref{controls} is represented by~\eqref{Volterra_epsilon} as follows:
$$
x_1(\varepsilon)= x_1^0 + \varepsilon v_1,\; x_2(\varepsilon)= x_2^0 + \varepsilon v_2,
$$
\begin{equation}
x_3(\varepsilon)= x_3^0 + \varepsilon \left(v_1 x_2^0 - v_2 x_1^0\right) - \frac{\varepsilon^2}{2\pi k_{12}} a_{12} (a_{12}-2v_1).
\label{x_representation}
\end{equation}
Note that representation~\eqref{x_representation} is exact (i.e. the higher order terms $R(\varepsilon)$ in the Volterra expansion vanish) as system~\eqref{Brockett3} is nilpotent.
This implies the following lemma.

\begin{lemma}
For arbitrary $x^0=(x_1^0,x_2^0,x_3^0)'\in {\mathbb R}^3$, $x^1=(x_1^1,x_2^1,x_3^1)'\in {\mathbb R}^3$, and $\varepsilon>0$, define the controls $u_1=u_1(t)$ and $u_2=u_2(t)$ by formulas~\eqref{controls} with
$$
v_1 = \frac{x_1^1-x_1^0}{\varepsilon},\; v_2 = \frac{x_2^1-x_2^0}{\varepsilon},
$$
$$
a_{12} = \frac{x_1^1-x_1^0}{\varepsilon} \pm \frac{1}{\varepsilon}\sqrt{(x_1^1-x_1^0)^2+2\pi k_{12} (x_3^0-x_3^1+x_1^1 x_2^0 -x_1^0 x_2^1)}.
$$
 Then the corresponding solution $x(t)$ of system~\eqref{Brockett3} with initial data $x(0)=x^0$ satisfies the condition $x(\varepsilon)=x^1$.
\end{lemma}

To solve the stabilization problem for system~\eqref{Brockett3}, consider a Lyapunov function candidate
$$
V(x) = \frac{1}{2}(x_1^2+x_2^2+x_3^2).
$$
Following the approach of Theorem~\ref{thm_sampling}, we define a time-varying feedback control to approximate the gradient flow of $-\nabla V(x)$ by trajectories of system~\eqref{Brockett3}:
\begin{equation}
u_1^\varepsilon (t,x) = v_1(x) + a(x)\cos\left(\frac{2\pi k (x)}{\varepsilon}t\right),
\label{u1_control}
\end{equation}
\begin{equation}
u_2^\varepsilon (t,x) = v_2(x) + a(x)\sin\left(\frac{2\pi k (x)}{\varepsilon}t\right),
\label{u2_control}
\end{equation}
where
$$
v_1(x) = -x_1,\; v_2(x) = -x_2,\; k(x) = {\rm sign}\,x_3,
$$
$$
a(x)=\left\{\begin{array}{ll}-x_1\pm \sqrt{x_1^2+\frac{2\pi |x_3|}{\varepsilon}},& x_3\neq 0, \\ 0, & x_3=0.\end{array}\right.
$$
 Without loss of generality, we may assume any integer value for $k(x)$ if $x_3=0$.

By Theorem~\ref{thm_sampling}, the feedback control~\eqref{u1_control}--\eqref{u2_control} ensures exponential stability of the equilibrium $x=0$ in the sense of $\pi_\varepsilon$-solutions,
provided that $\varepsilon>0$ is small enough.

\section{Simulation results}

\begin{figure}
\includegraphics[angle=270]{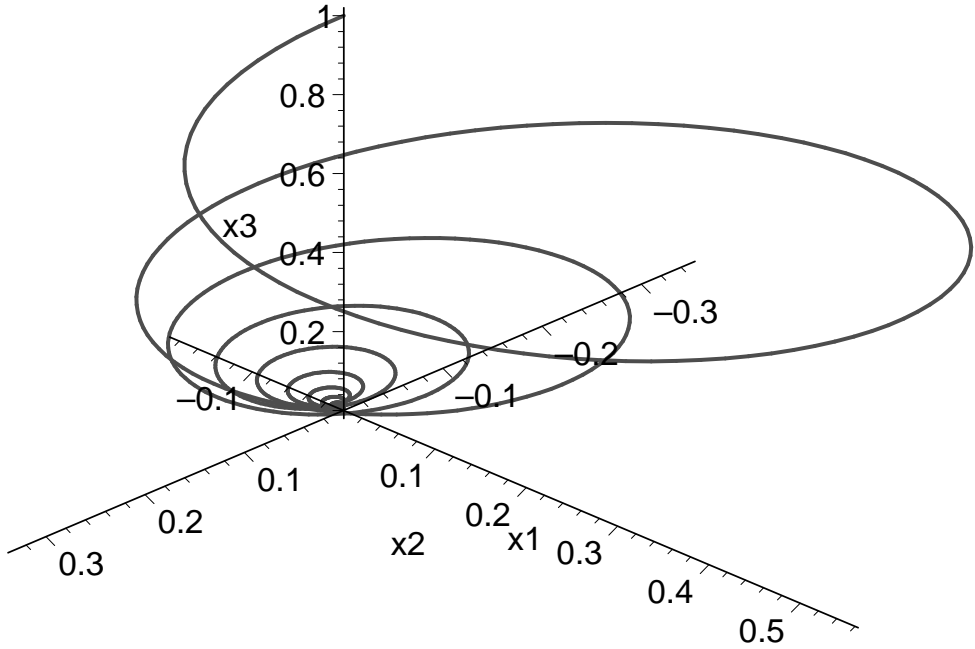}
\caption{Trajectory of the closed-loop system~\eqref{Brockett3}, \eqref{u1_control}, \eqref{u2_control} for $x_1(0)=x_2(0)=0$, $x_3(0)=1$.}
\label{fig:1}
\end{figure}

\begin{figure}
\includegraphics[scale=0.75,angle=270]{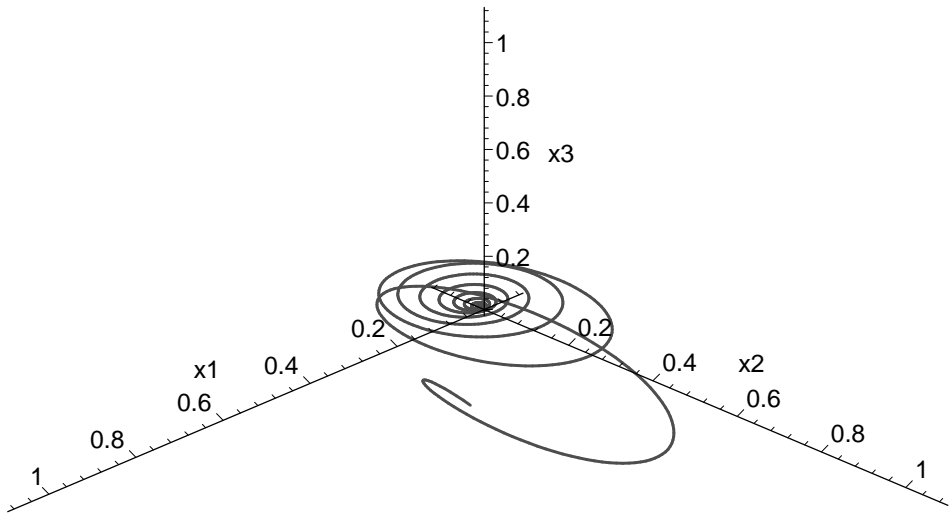}
\caption{Trajectory of the closed-loop system~\eqref{Brockett3}, \eqref{u1_control}, \eqref{u2_control} for $x_1(0)=x_2(0)=x_3(0)=1$.}
\label{fig:2}
\end{figure}

In this section,
we perform numerical integration of the closed-loop system~\eqref{Brockett3}
with the feedback law $u=u(t,x(t))$ of form~\eqref{u1_control}--\eqref{u2_control}.
Trajectories of this system are shown in Fig.~\ref{fig:1} and~\ref{fig:2} for $\varepsilon=1$ and the following initial conditions:
$$
x_1(0) = x_2(0) = 0,\; x_3(0) = 1\quad \text{(Fig.~\ref{fig:1})},
$$
$$
x_1(0) = x_2(0) = x_3(0) = 1\quad \text{(Fig.~\ref{fig:2})}.
$$


These simulation results show that the feedback law~\eqref{u1_control}--\eqref{u2_control} steers the Brockett integrator to the origin not only in the sense of $\pi_\varepsilon$-solutions (as stated in Theorem~\ref{thm_sampling}), but also in the sense of classical solutions.

\section{Conclusion}

In this paper, a family of time-dependent trigonometric polynomials with coefficients depending on the state has been constructed
to stabilize the equilibrium of a nonholonomic system.
These coefficients are obtained by solving an auxiliary system of quadratic algebraic equations involving the gradient of a Lyapunov function.
An important feature of this work relies on the proof of the solvability of such a system for an arbitrary dimension of the state space provided that the Lie algebra rank condition is satisfied with first order Lie brackets.
It should be emphasized that this result is heavily based on the degree principle as the implicit function theorem is not applicable for a non-differentiable function $G_x(\xi)$ in Lemma~\ref{lemma_solvability}.

Another important remark is that our design scheme produces small controls $u^\varepsilon(t,x)$ for small values of $\|x\|$, and the frequencies of the sine and cosine functions are constant for each fixed $\varepsilon>0$.
This feature differs from the approach to the motion planning problem that uses a sequence of high-amplitude highly oscillating open-loop controls (see~\cite{SL,Liu}).

The proof of Theorem~\ref{thm_sampling} is considered as an extension of Lyapunov's direct method, where the decay condition for a Lyapunov function is guaranteed by exploiting the Volterra expansion instead of using the time derivative along the trajectories. Although the exponential stability result is established for $\pi_\varepsilon$-solutions under a sampling strategy, simulation results demonstrate the convergence of classical solutions of the closed-loop system to its equilibrium.
Thus, the question of the limit behavior of classical (or Carath\'eodory) solutions of system~\eqref{Sigma} with the feedback control~\eqref{S1} remains open for further theoretical studies.

\section*{Acknowledgement}

The author is grateful to Prof. Bronis{\l}aw Jakubczyk for valuable discussions.

\end{document}